\documentclass[twoside]{amsart}

\usepackage{amssymb}
\usepackage{amsmath}
\usepackage{cite}
\usepackage{mathrsfs}

\usepackage{pstricks,pst-node,pst-coil}
\psset{linewidth=0.4pt}  

\usepackage{booktabs}
\makeatletter

\address{\bigskip\hfil\begin{tabular}{l@{}}
            School of Mathematics and Statistics F07\\  
            University of Sydney,
            Sydney N.S.W. 2006.\hfill\qquad
                        {\tt mathas@maths.usyd.edu.au}\\  
            Australia.\hfill\qquad
                    {\tt www.maths.usyd.edu.au/u/mathas/}
          \end{tabular}}

\copyrightinfo{}{}

\renewcommand{\labelenumi}{(\roman{enumi})}

\let\atop\@@atop

\def\And{\text{\ and\ }}

\def\ForAll{\text{\ for all\ }}

\def\If{\text{\ if\ }}

\def\th{{\text{th}}}

\def\Set[#1]#2|#3|{\Big\{\ #2\ \Big| \
            \vcenter{\hsize #1mm\centering#3}\Big\}}

{\catcode`\|=\active
  \gdef\set#1{\mathinner{\lbrace\,{\mathcode`\|"8000%
                                   \let|\midvert #1}\,\rbrace}}
}
\def\midvert{\egroup\mid\bgroup}

\def\Number#1{\refstepcounter{equation}
              \leqno(\theequation)\if*#1%
              \else\def\@currentlabel{{\rm\theequation}}\label{#1}%
              \fi}
\def\Dag{\ifmmode\leqno(\dag)\else$(\dag\)$\fi}
\def\DDag{\ifmmode\leqno(\ddag)\else$(\ddag\)$\fi}

\numberwithin{equation}{section}

\swapnumbers

\newtheorem{Theorem}[equation]{Theorem}
\newtheorem{Proposition}[equation]{Proposition}
\newtheorem{Lemma}[equation]{Lemma}
\newtheorem{Corollary}[equation]{Corollary}

\theoremstyle{remark}
\newtheorem{Remark}[equation]{Remark}

\newtheorem{Example}[equation]{Example}

  {\ifx*#1\let\pointlabel\relax\else\def\pointlabel{#1}\fi
   \refstepcounter{equation}\trivlist
   \item[\hskip\labelsep\bf\theequation
         \ifx\pointlabel\relax\else\space\pointlabel\space\fi]
   \ignorespaces\it
  }{\relax}

\let\MC\multicolumn

\def\){\big)}
\def\({\big(}
\let\>\rangle
\let\<\langle

\let\ss\subseteq
\let\0\varnothing              

\let\realb@r\bar
\let\bar\overline

\let\gedom\trianglerighteq
\let\gdom\vartriangleright

\def\F{{\mathbb F}}
\def\N{{\mathbb N}}
\def\Z{{\mathbb Z}}

\def\C{{\mathbb C}}

\def\map#1#2{\,{:}\,#1\!\longrightarrow\!#2}

\author{Gordon James and Andrew Mathas} 
\title{Equating decomposition numbers for different~primes} 
\subjclass[2000]{20C08, 20C30, 33D80}
\def\addresses{
\bigskip\makeatletter\sc
  \begin{tabular}[t]{ll}
  G.~J.& Department of Mathematics, Imperial College\\
       & 180 Queen's Gate, London SW7 2BZ. U.K.\\
       & {\tt g.james@ic.ac.uk}\hfill
         {\tt www.ma.ic.ac.uk/$\sim$gdj/}\\[5pt]
  A.~M.&School of Mathematics F07, University of Sydney\\
       &Sydney N.S.W. 2006.\hfill\qquad
	                {\tt a.mathas@maths.usyd.edu.au}\\
       &Australia.\hfill\qquad
	            {\tt www.maths.usyd.edu.au/u/mathas/}
\end{tabular}
}

\def\A{\mathcal A}

\def\aF{{}^\alpha\!F}
\def\F{\mathscr F}
\def\H{\mathscr H}
\def\Sc{\mathscr S}
\def\Hei{\mathcal H_e}
\def\Hm{\mathcal H^-_e}
\def\Hmp{\mathcal H^-_{e+1}}

\def\len(#1){\ell(#1)}
\def\Rad{\operatorname{Rad}}
\def\res{\operatorname{res}}
\def\Sym{\mathfrak S}
\def\Ue{{\bf U}_v(\widehat{\mathfrak{sl}}_e)}
\def\UeA{{\bf U}_\A(\widehat{\mathfrak{sl}}_e)}
\def\Uep{{\bf U}_v(\widehat{\mathfrak{sl}}_{e+1})}
\def\Ume{{\bf U}^-_v(\widehat{\mathfrak{sl}}_e)}
\def\UmeA{{\bf U}^-_\A(\widehat{\mathfrak{sl}}_e)}
\def\Umep{{\bf U}^-_v(\widehat{\mathfrak{sl}}_{e+1})}
\def\Upe{{\bf U}^+_v(\widehat{\mathfrak{sl}}_e)}

\def\V{\mathcal V}

\def\head{\operatorname{head}}
\def\spin{\operatorname{spin}}
\def\ribbon#1{\,\everypsbox{\scriptstyle}
         \pccoil[linewidth=0.3pt,coilarm=.17,coilheight=1.2,coilwidth=.1]
                {->}(0,0.1)(1,0.1)
         \rput(.5,.25){#1}
         \hspace*{11mm}\vrule height 2.2ex width 0pt
}
\def\sribbon#1{\,\everypsbox{\scriptscriptstyle}
         \pccoil[coilarm=.15,coilheight=1.2,coilwidth=.08]
                {->}(0,0.07)(0.8,0.07)
         \rput(0.4,.19){#1}
         \hspace*{8mm}\vrule height 1.6ex width 0pt
}

\def\myarrow#1{\,\ooalign{\raise1.2ex\hbox{\scriptsize\ \ $#1$}\crcr%
           $\xrightarrow{\phantom{\hbox{\scriptsize\ \ $#1$}}}$}\,}
\def\iarrow{\myarrow{i}}

\let\section=\subsection
\def\theequation{\arabic{subsection}.\arabic{equation}}
\def\thesection#1{\relax}
\makeatletter
\@addtoreset{equation}{subsection}
\makeatother

\begin{document}
\begin{abstract}
This paper shows that certain decomposition numbers for the
Iwahori--Hecke algebras of the symmetric groups and the $q$--Schur
algebras at different roots of unity in characteristic zero are equal.
To prove our results we  first establish the corresponding theorem for
the canonical basis of the level one Fock space and then apply deep
results of Ariki and Varagnolo and Vasserot.
\end{abstract}
 
\maketitle

\section{Introduction}
Throughout this note we adopt the standard notation for the modular
representation theory of the symmetric groups, as can be found in
\cite{James,M:ULect}.

Consider the following two submatrices of the $p$--modular decomposition
matrices of the symmetric groups $\Sym_n$.
$$\begin{array}{l|*{9}{l}}
10,1   &1&&&&&&&&\\
9,2    &1&1&&&&&&&\\
7,4    &.&1&1&&&&&&\\
7,2^2  &1&1&1&1&&&&&\\
6,5    &.&.&1&.&1&&&&\\
6,2^2,1&1&1&1&1&1&1&&&\\
4^2,3  &.&.&1&1&1&.&1&&\\
4^2,2,1&.&1&1&1&\bf2&1&1&1\\
\multicolumn{10}c{n=11\And p=3}\\
\end{array}
\quad
\begin{array}{l|*{9}{l}@{}}
18,3    &1&&&&&&&&\\
17,4    &1&1&&&&&&&\\
13,8    &.&1&1&&&&&&\\
13,4^2  &1&1&1&1&&&&&\\
12,9    &.&.&1&.&1&&&&\\
12,4^2,1&1&1&1&1&1&1&&&\\
8^2,5   &.&.&1&1&1&.&1&&\\
8^2,4,1 &.&1&1&1&\bf1&1&1&1\\
\multicolumn{10}c{n=21\And p=5}\\
\end{array}$$
The two matrices are identical except for the labelling and the two bold 
faced entries (omitted entries are zero).

This paper was motivated by our attempts to compute the bold faced
entry in the matrix above for $n=21$ and $p=5$. At the outset we knew
that this number was either $1$ or $2$; however, we were unable to
determine which of these possibilities was correct. L\"ubeck and
M\"uller~\cite{LM:p=5} have shown that this multiplicity is
equal to~$1$ using computer calculations; see~\cite[\S5.3]{LN:p=5} 
for details.

Note that the partitions $(4^2,2,1)$ and $(8^2,4,1)$ both have
$p$--weight $3$ and $p$--core $(p-2,p-2)$ for $p=3$ and $p=5$
respectively. Martin and Russell~\cite{MartinRussell} have claimed that
all of the $p$--modular decomposition numbers of the symmetric group
of $p$--weight~$3$ are~$0$ or~$1$ when $p>3$; unfortunately, their
proof contains a gap when dealing with partitions with $p$--core
$(p-2,p-2)$. This particular case is still open when $p>5$; as a
consequence, the claim in~\cite{MartinRussell} that the decomposition
numbers are always $0$ or $1$ for partitions of $p$--weight~$3$ when
$p>3$ is in doubt.

In this paper we prove a general theorem which indicates why the
matrices above are very similar. This result is not about the
decomposition matrices of the symmetric groups but rather about the
closely related decomposition matrices of the $q$--Schur algebras
$\Sc_{\C,q}(n)$ at a complex root of unity. Our result shows that
certain decomposition numbers of $\Sc_{\C,q}(n)$ and $\Sc_{\C,q'}(m)$
are equal for specified $m>n$. The decomposition matrix for the
Iwahori--Hecke algebra $\H_{\C,q}(\Sym_n)$ is a submatrix of the
decomposition matrix of $\Sc_{\C,q}(n)$; so, in particular, our result
shows that for the Iwahori--Hecke algebras all of the decomposition
multiplicities above for $(n,e)=(11,3)$ and $(m,e')=(21,5)$ are equal
(in the Hecke algebra case the multiplicities $d_{(4^2,2,1),(6,5)}$,
when $e=3$, and $d_{(8^2,4,1),(12,9)}$, when $e=5$, are both equal
to~$1$).

\section{Abacuses and the $q$--Schur algebra} In order to state
our results we recall the abacus notation for partitions introduced
in~\cite{James:YoungD}. Fix an integer $e\ge2$. An {\sf $e$--abacus}
is an abacus with $e$ runners, which we label from left to right as
$\rho_0,\dots,\rho_{e-1}$. Number the bead positions on the abacus by
$0,1,2,\dots$ reading from left to right and then top to bottom; so
the bead positions on $\rho_r$ are numbered $r+me$ for $m\ge0$. We
order the beads on a given $e$--abacus according to their bead
positions.

Let $\lambda=(\lambda_1,\lambda_2,\dots)$ be a partition. The 
{\sf length} of $\lambda$ is the smallest integer $\len(\lambda)$ such
that $\lambda_i=0$ for all $i>\len(\lambda)$. If $k\ge\len(\lambda)$
then $\lambda$ has a (unique) $e$--abacus representation with $k$
beads; namely, the $e$--abacus with beads at positions $$\lambda_k,
\lambda_{k-1}+1,\dots,\lambda_2+k-2,\lambda_1+k-1.$$ The bead
positions on an abacus for $\lambda$ encode the first column hook
lengths, so this gives a natural bijection between the abacuses with
$k$ beads and the partitions of length at most $k$.  For our purposes
it is important that the $m^\th$ bead on an abacus for $\lambda$
corresponds to row $l=k-m+1$ of $\lambda$ (row $l$ of $\lambda$ is
empty if $l>\len(\lambda)$). 

In what follows we fix $k\ge0$ and (with few exceptions) consider
only abacuses with $k$ beads; or, equivalently, partitions of length
at most $k$.

We want to compare $e$--abacuses with $(e+1)$--abacuses. Fix an
integer $\alpha$ with $0\le\alpha\le e$. If $\lambda$ is a partition
with $\len(\lambda)\le k$ then $\lambda$ can be represented on an
$e$--abacus with $k$ beads. Let $\lambda^+$ be the partition
corresponding to the $(e+1)$--abacus obtained by inserting an {\it
empty} runner before $\rho_\alpha$ in the $e$--abacus for $\lambda$
(if $\alpha=e$ we insert an empty $e^\th$ runner). Let
$\rho_0^+,\dots,\rho_e^+$ be the runners of the $(e+1)$--abacus of
$\lambda^+$; then $\rho^+_r=\rho_r$ if $r<\alpha$, $\rho^+_\alpha$
is empty, and $\rho^+_r=\rho_{r-1}$ if $r>\alpha$.

Although our notation does not reflect this the partition $\lambda^+$
does depend upon both the choice of $\alpha$ and the choice of $k$. 

\begin{Example}
Suppose that $e=3$, $k=4$ and $\alpha=2$. Let
$\lambda=(4^2,3)$. Then the abacuses (with $4$ beads) for $\lambda$,
$\lambda^+$ and $\lambda^{++}=(\lambda^+)^+$ are as follows.
$$\begin{array}{@{\hskip 14mm}*3{c@{\hskip 12mm}}}
\begin{array}[t]{*3l}\MC3c{e=3}\\
  2 & 0 &1\\\hline
  \bullet&\cdot&\cdot\\
  \cdot&\bullet&\cdot\\
  \bullet&\bullet&\cdot\\
  \cdot&\cdot&\cdot
\end{array} & 
\begin{array}[t]{*4c}\MC4c{e'=4}\\
  0&1&2&3\\\hline
  \bullet&\cdot&\cdot&\cdot\\
  \cdot&\bullet&\cdot&\cdot\\
  \bullet&\bullet&\cdot&\cdot\\
  \cdot&\cdot&\cdot&\cdot
\end{array} &
\begin{array}[t]{*5c}\MC5c{e''=5}\\
  1&2&3&4&0\\\hline
  \bullet&\cdot&\cdot&\cdot&\cdot\\
  \cdot&\bullet&\cdot&\cdot&\cdot\\
  \bullet&\bullet&\cdot&\cdot&\cdot\\
  \cdot&\cdot&\cdot&\cdot&\cdot
\end{array} \\[8pt]
\lambda=(4^2,3) & \lambda^+=(6^2,4) & \lambda^{++}=(8^2,5)
\end{array}$$
We have labelled the runners by their residues which will be introduced
below.

The reader is invited to check that the partitions which label the
decomposition matrix for $n=21$ and $p=5$ in the introduction are
precisely the partitions $\lambda^{++}$ as $\lambda$ runs over the
corresponding partitions of $11$. We emphasize that the empty runner
can be inserted anywhere in the abacus.
\label{example}\end{Example}

We are now almost ready to describe our main result. As in the
introduction let~$q$ be a primitive $e^\th$ root of unity in $\C$ and
let $\Sc_{\C,q}(n)$ be the $q$--Schur algebra defined over the complex
numbers with parameter $q$; so $\Sc_{\C,q}(n)=\Sc_{\C,q}(n,n)$ in the
notation of Dipper and James~\cite{DJ:Schur}.

For each partition $\lambda$ of $n$ Dipper and James~\cite{DJ:qWeyl}
(see also \cite{M:ULect}), defined a right $\Sc_{\C,q}(n)$--module
$W^\lambda_q$, called a {\sf Weyl module}. There is a natural bilinear
form $\<\ ,\ \>$ on $W^\lambda_q$ and 
$\Rad W^\lambda_q
   =\set{x\in W^\lambda_q|\<x,y\>=0\ForAll y\in W^\lambda_q}$
is an $\Sc_{\C,q}(n)$--submodule of~$W^\lambda_q$; set
$L^\lambda_q=W^\lambda_q/\Rad W^\lambda_q$.  Dipper and James showed
that $L^\lambda_q$ is an absolutely irreducible $\Sc_{\C,q}(n)$--module
and, further, that every irreducible $\Sc_{\C,q}(n)$--module arises
uniquely in this way.  Let $[W^\lambda_q:L^\mu_q]$ be the multiplicity
of the simple module $L^\mu_q$ as a composition factor of
$W^\lambda_q$.

Let $q'$ be a primitive $(e+1)^{\text{st}}$ root of unity in $\C$.
Then we also have the $q'$--Schur algebra $\Sc_{\C,q'}(m)$ and modules
$W^\nu_{q'}$ and $L^\nu_{q'}$, for $\nu$ a partition of $m$. We shall
prove the following.

\begin{Theorem}
Suppose that $\lambda$ and $\mu$ are partitions of $n$ of
length at most $k$. Then 
$$[W^\lambda_q:L^\mu_q]=[W^{\lambda^+}_{q'}:L^{\mu^+}_{q'}].$$
\label{main}\end{Theorem}

It may not be clear to the reader that this result is really saying
that the decomposition matrices of the blocks containing $W^\lambda_q$
and $W^{\lambda^+}_{q'}$ are equal on those rows indexed by partitions
of length at most $k$, when we order the rows of these matrices in a
way compatible with dominance. This follows from
Lemma~\ref{bookkeeping} below.

Let $\H_{\C,q}(\Sym_n)$ be the Iwahori--Hecke algebra of
$\Sym_n$~\cite{M:ULect}. Then for each partition $\lambda$ there is an
$\H_{\C,q}(\Sym_n)$--module $S^\lambda_q$, called a Specht module,
which carries an associative bilinear form. Let
$D^\lambda_q=S^\lambda_q/\Rad S^\lambda_q$; then $D^\lambda_q$ is
either zero or absolutely irreducible and every irreducible
$\H_{\C,q}(\Sym_n)$--module arises uniquely in this way. Moreover,
$D^\lambda_q\ne0$ if and only if $\lambda$ is $e$--regular; that is,
if and only if no $e$ non--zero parts of $\lambda$ are equal.

There is a $q$--analogue of the Schur functor which maps $W^\lambda_q$
to $S^\lambda_q$ and $L^\lambda_q$ to $D^\lambda_q$ for each
$\lambda$; in particular, this shows that
$[W^\lambda_q:L^\mu_q]=[S^\lambda_q:D^\mu_q]$ whenever $\mu$ is
$e$--regular. Hence, from Theorem~\ref{main} we obtain the following:

\begin{Corollary}
Suppose that $\lambda$ and $\mu$ are partitions of $n$ of
length at most $k$ such that $\mu$ is $e$--regular. Then 
$[S^\lambda_q:D^\mu_q]=[S^{\lambda^+}_{q'}:D^{\mu^+}_{q'}].$
\label{regular main}\end{Corollary}

It is tempting to speculate that there is some form of category
equivalence underpinning these results. However, in general, there are
a different number of simple modules in the blocks for $\lambda$ and
$\lambda^+$, so these blocks are certainly not Morita equivalent.

Rather than prove our comparison theorem for decomposition numbers
directly we prove a stronger result relating the
LLT--polynomials~\cite{LLT,LT:Schur2}. We now recall the notation
needed to describe this.

\section{The Fock space and $\Ue$ --- the regular case}
Let $v$ be an indeterminate over $\C$. The {\sf Fock space} is the infinite
rank free $\C[v,v^{-1}]$--module
$\F=\bigoplus_{n\ge0}\bigoplus_{\lambda\vdash n}\C[v,v^{-1}]\lambda$.
The Fock space has a natural structure as a module for the affine
quantum group $\Ue$; we will describe how the negative part
$\Ume$ of $\Ue$ acts on $\F$ since this is all we shall need (full
details can be found in \cite{LLT,M:ULect}). 

The {\sf diagram} of a partition $\lambda$ is the set
$[\lambda]=\set{(c,d)\in\N^2|d\le\lambda_c}$. A {\sf node} is any ordered
pair $(c,d)\in\N^2$; in particular, all of the elements of $[\lambda]$
are nodes.  The {\sf $e$--residue} of the node $x=(c,d)$ is
$\res_e(x)=d-c\pmod e$; $x$ is an $i$--node if
$res_e(x)=i$.

A node $x$ is an {\sf addable} node of $\lambda$ if $[\lambda]\cup\{x\}$
is the diagram of a partition (and~$x\notin[\lambda]$); similarly, $x$
is {\sf removable} if $[\lambda]\setminus\{x\}$ is the diagram of a
partition (and~$x\in[\lambda]$). For $i=0,\dots,e-1$ let
$A_i(\lambda)$ be the set of addable $i$--nodes for~$\lambda$ 
and~$R_i(\lambda)$ be the set of removable $i$--nodes.  Given two nodes
$x=(c,d)$ and~$y=(a,b)$ say that $y$ is {\sf above} $x$ if $c>a$; if $y$
is above $x$ we write $y\succ x$.

In order to define the action of $\Ume$ on $\F$ for
$i=0,\dots e-1$ write $\lambda\iarrow\nu$ if $\nu$ is a partition of
$n+1$ and $[\nu]=[\lambda]\cup\{x\}$ for some addable $i$--node $x$.
Finally, we set
$N_i(\lambda,\nu)=\#A_i(\lambda,\nu)-\#R_i(\lambda,\nu)$ where
$A_i(\lambda,\nu)=\set{y\in A_i(\lambda)|y\succ x}$ and
$R_i(\lambda,\nu)=\set{y\in R_i(\lambda)|y\succ x}$. 

Let $F_0,\dots,F_{e-1}$ be the Chevalley generators of $\Ume$. Then
the action of~$\Ume$ on $\F$ is determined by 
$$ F_i\lambda=\sum_{\lambda\iarrow\nu}v^{N_i(\lambda,\nu)}\nu,
\Number{U-action}$$ 
for $0,1,\dots,e-1$.

Let $\Lambda_0,\dots,\Lambda_{e-1}$ be the fundamental weights of
$\Ue$ and let $L(\Lambda_0)$ be the
irreducible integrable highest weight module of high weight $\Lambda_0$. 
Then $L(\Lambda_0)\cong\Ue\0=\Ume\0$ as a $\Ue$--module~\cite{LLT},
where $\0\in\F$ is the empty partition.

Let $\rule[1.4ex]{.6em}{.1ex}$\ be the bar involution on $\UeA$,
Kostant--Lusztig $\A$--form of $\Ue$ (where $\A=\Z[v,v^{-1}]$). Then
Lascoux, Leclerc and Thibon~\cite{LLT} showed that $L(\Lambda_0)$ has
a basis $\set{B_\mu|\mu\ e\text{--regular}}$ which is uniquely
determined by the requirements that $\bar B_\mu=B_\mu$ and
$$B_\mu=\sum_{\substack{\lambda\vdash n\\\mu\gedom\lambda}}
                 b_{\lambda\mu}(v)\lambda$$
for some polynomials $b_{\lambda\mu}(v)\in\Z[v]$ such that
$b_{\mu\mu}(v)=1$ and $b_{\lambda\mu}(v)\in v\Z[v]$ whenever
$\lambda\ne\mu$. This basis is the Kashiwara--Lusztig canonical basis
of $L(\Lambda_0)$.

In particular, note that Lascoux, Leclerc and Thibon~\cite{LLT}
showed that $b_{\lambda\mu}(v)=0$ if either $|\lambda|\ne|\mu|$ or if
$\lambda$ and $\mu$ have different $e$--cores.

We want to compare the actions of $\Ue$ and $\Uep$ on the Fock space
$\F$. In order to distinguish between these two algebras let
$F_0^+,\dots,F_e^+$ be the Chevalley generators of $\Umep$, let
$\Lambda_0^+,\dots,\Lambda_e^+$ be its fundamental weights, and let
$\F^+\cong\F$ be the Fock space for the $\Uep$--action. Then
$L(\Lambda^+_0)\cong\Umep\0$ as a $\Uep$--module. Given an
$(e+1)$--regular partition $\nu$ let
$$B^+_\nu=\sum_{\nu\gedom\sigma} b^+_{\sigma\nu}(v)\sigma$$
be the corresponding canonical basis element of $L(\Lambda^+_0)\ss\F^+$.

We can now state a stronger version of Corollary~\ref{regular main}.

\begin{Theorem}
Suppose that $\lambda$ and $\mu$ are partitions
of $n$ of length at most $k$ and suppose that
$\mu$ is $e$--regular. Then
$b_{\lambda\mu}(v)=b_{\lambda^+\mu^+}^+(v)$.
\label{fock}\end{Theorem}

Ariki \cite[Prop.~4.3(2)]{Ariki:can} has shown that the polynomials
$b_{\lambda\mu}(v)$ at $v=1$ compute the decomposition multiplicities;
explicitly, $[S^\lambda_q:D^\mu_q]=b_{\lambda\mu}(1)$ and
$[S^{\lambda^+}_{q'}:D^{\mu^+}_{q'}]=b_{\lambda^+\mu^+}^+(1)$.
Consequently, Theorem~\ref{fock} implies Corollary~\ref{regular main}.
The result also hints at additional structure because, conjecturally,
the polynomials $b_{\lambda\mu}(v)$ and $b_{\lambda^+\mu^+}^+(v)$ also
describe the Jantzen filtrations of $S^\lambda_q$ and
$S^{\lambda^+}_{q'}$; see \cite{LLT,JM:Schaper}.

We prove Theorem~\ref{fock} directly using the LLT algorithm; in the
next section we will extend this argument to cope with the case where
$\mu$ is not necessarily $e$--regular.  

Fred Goodman has pointed out that Theorem~\ref{fock} can also be
deduced from~\cite[Theorem~5.3]{GW}. We remark that the origin of our
results, and those of Goodman and Wenzl, is that the
$b_{\lambda\mu}(v)$ are parabolic Kazhdan--Lusztig polynomials for the
parabolic subgroup $\Sym_k$ of the extended affine Weyl group
$\hat\Sym_k$ \cite{GW,VV,LT:Schur2}; in turn, the parabolic
Kazhdan--Lusztig polynomials are naturally indexed by the alcoves and,
generically, the alcove geometry does not depend on $k$ or $e$.

We begin the proof of Theorem~\ref{fock} with the following Lemma
which is largely book keeping. For example, the result implicitly
assumes that $\mu^+$ is $(e+1)$--regular.

\begin{Lemma}
Let $\lambda$ and $\mu$ be partitions of $n$ of length at
most $k$. Then
\begin{enumerate}
\item $\mu^+$ is $(e+1)$--regular;
\item $\lambda$ and $\mu$ have the same $e$--core if and only if
$\lambda^+$ and $\mu^+$ have the same $(e+1)$--core;
\item if $\lambda$ and $\mu$ have the same $e$--core then
$|\lambda^+|=|\mu^+|$; and,
\item $\lambda\gedom\mu$ if and only if $\lambda^+\gedom\mu^+$
\end{enumerate}\label{bookkeeping}
\end{Lemma}

\begin{proof}
A partition is $(e+1)$--regular if and only if its
$(e+1)$--abacus does not contain a string of $e+1$
consecutive beads. Hence, $\mu^+$ is $(e+1)$--regular since the runner
$\rho_\alpha^+$ is empty; this proves (i). (In fact, if $\mu$ is
$e$--regular then so is~$\mu^+$.)

Next, recall that the $e$--abacus for the $e$--core of $\lambda$ is
obtained by rearranging the beads on each runner of the $e$--abacus
for $\lambda$ in such a way that no bead has an empty bead position
above it. Hence, if $\kappa$ is the $e$--core of $\lambda$ then
$\kappa^+$ is the $(e+1)$--core of~$\lambda^+$, so (ii) follows.

For (iii) define $w$ by $|\lambda|=|\kappa|+we$; in
other words, $w$ is the $e$--weight of~$\lambda$. Now, $w$ can be read
off the $e$--abacus for $\lambda$ by adding up, for each bead $\beta$,
the number of empty bead positions which are above~$\beta$ and also
on the same runner.  Consequently, $w$ is also the $(e+1)$--weight of
$\lambda^+$; hence, $|\lambda^+|=|\kappa^+|+w(e+1)$. Lastly, since
$\kappa$ is also the $e$--core of $\mu$ it follows that $\mu$ is also
a partition of $e$--weight~$w$ and that
$|\mu^+|=|\kappa^+|+w(e+1)=|\lambda^+|$, as required.

Finally, let $\beta_1^\lambda<\dots<\beta^\lambda_k$
and  $\beta_1^\mu<\dots<\beta^\mu_k$ be the positions of the beads on
the $e$--abacuses for $\lambda$ and $\mu$ respectively. Then it is
easy to see that $\lambda\gedom\mu$ if and only if 
$\sum_{s=t}^k\beta^\lambda_s\ge\sum_{s=t}^k\beta^\mu_s$
for $t=1,\dots,k$. Rephrasing this condition in terms of the
$(e+1)$--abacuses of $\lambda^+$ and $\mu^+$ proves (iv).
\end{proof}

We remark that if $\lambda$ and $\mu$ are partitions of $n$ with
different $e$--cores then, in general, it is {\it not} true that
$|\lambda^+|=|\mu^+|$. 

Let $\F_{>k}$ be the $\C[v,v^{-1}]$--submodule of $\F$ spanned by the
partitions of length strictly greater than $k$. By (\ref{U-action})
$\F_{>k}$ is a $\Ume$--submodule of $\F$; it is not, however, a
$\Ue$--submodule. Therefore, $\F_k=\F/\F_{>k}$ is a $\Ume$--module.
We abuse notation and identify the elements of $\F$ with their images
in $\F_k$; with this understanding, $\set{\lambda|\len(\lambda)\le k}$
is a basis of $\F_k$. 

Similarly, $\F_k^+=\F^+/\F_{>k}^+$ is a $\Umep$--module.  We want to
compare the action of $\Ume$ on $\F_k$ with the action of $\Umep$ on
$\F^+_k$; to do this we reinterpret (\ref{U-action}) in terms of
abacuses.

Suppose $\lambda$ is a partition with $\len(\lambda)\le k$ and
consider the $e$--abacus of $\lambda$. For $0\le r<e$ define the
{\sf $e$--residue} of the runner $\rho_r$ to be the integer
$\res_e(\rho_r)$ determined by the following two conditions.
\begin{enumerate}
\item The $e$--residue of the runner which holds the last bead is
$\lambda_1-1\pmod e$.
\item Modulo $e$, the $e$--residues of the runners increase by $1$ from
left to right.
\end{enumerate}\noindent%
In Example~\ref{example} the runners are labelled by their $f$--residues for
$f=3,4$ and~$5$ respectively. Similarly, we define the
$(e+1)$--residues $\res_{e+1}(\rho^+_r)$, for $0\le r\le e$, of the
runners of the $(e+1)$--abacus for $\lambda^+$.

The $e$--residue of a bead $\beta$ is defined to be the $e$--residue
of the corresponding runner. As we have seen, the $k$ beads on the
$e$--abacus correspond to the first $k$ rows of $\lambda$ (in reverse
order); it is easy to see that the $e$--residue of a bead is equal to
the $e$--residue of the node at the end of the corresponding row of
$\lambda$. In particular, this implies that $e$--residues of the
runners depend only on the $e$--core of $\lambda$.

The operator $F_i$ acts on a partition $\lambda$ by adding nodes of
$e$--residue~$i$. Because the $e$--residues of the runners correspond
to the $e$--residues of nodes at the end of the rows of $\lambda$,
this is the same as moving a bead on the $e$--abacus of $\lambda$ from
the runner with $e$--residue $i-1$ to an adjacent empty position on
the runner with $e$--residue $i$.

Recall that in the definition of $\lambda^+$ we have fixed an integer
$\alpha $ with $0\le \alpha<e$. We now introduce a
$\C[v,v^{-1}]$--linear map $\aF_i\map{\F_k}{\F^+_k}$ for
$i=0,\dots,e-1$. To define $\aF_i$ it is enough to describe
$\aF_i\lambda$ for each partition~$\lambda$ with $\len(\lambda)\le k$.
As above, let $\rho_0,\dots,\rho_{e-1}$ be the runners of the
$e$--abacus of $\lambda$ (with~$k$ beads) and let
$\rho^+_0,\dots,\rho^+_e$ be the runners of the $(e+1)$--abacus for
$\lambda^+$. There is a unique $r$ such that $i=\res_e(\rho_r)$ (and
$0\le r<e$); set $j=\res_{e+1}(\rho^+_r)$. Define
$$\aF_i\lambda=\begin{cases} 
                   F^+_j\lambda^+,&\If 0\le r<\alpha,\\
                   F^+_{j+1}F^+_j\lambda^+,&\If r=\alpha,\\
                   F^+_{j+1}\lambda^+,&\If \alpha<r<e,
\end{cases}$$
where $j+1$ is understood modulo $e$. Similarly, we define the divided
powers $\aF_i^{(a)}$ for $a\ge1$; for example, when $r=\alpha$ we set
$\aF_i^{(a)}\lambda=F_{j+1}^{+(a)}F_j^{+(a)}\lambda^+$.

For our final piece of notation observe that
$\len(\lambda^+)\ge\len(\lambda)$ for any partition $\lambda$ (and if
$\lambda$ is a partition of length at most $k$ then
$\len(\lambda)\le\len(\lambda^+)\le k$). Therefore, we have a
well--defined $\C[v,v^{-1}]$--linear map $\Theta\map{\F_k}\F_k^+$
determined by $\Theta(\lambda)=\lambda^+$ for $\len(\lambda)\le k$. As
with $\lambda^+$, we emphasize that $\Theta$ depends upon
both~$\alpha$ and~$k$. The map~$\Theta$ is injective but not
surjective, having image the span of those partitions of length at
most $k$ which have an $(e+1)$--abacus with $k$ beads and with an
empty runner $\rho^+_\alpha$.

\begin{Lemma}
\label{commuting1}
Suppose that $0\le i<e$. Then the following diagram commutes.
$$\begin{array}{c}\\[-6pt]
\psmatrix[colsep=1cm,rowsep=1cm]
  \F_k & \F_k \\
       & \F^+_k
\endpsmatrix
\psset{nodesep=3pt,arrows=->}
\ncline{1,1}{1,2}\taput{F_i}
\ncline{1,1}{2,2}\tlput{\aF_i}
\ncline{1,2}{2,2}\trput{\Theta}
\end{array}$$
\end{Lemma}

\begin{proof}
It suffices to verify the lemma for a partition
$\lambda\in\F_k$. As above, let $\rho_r$ be the runner in the
$e$--abacus for $\lambda$ for which $i=\res_e(\rho_r)$ and set
$j=\res_{e+1}(\rho_r^+)$. 

First consider $F_i\lambda=\sum_\nu v^{N_i(\lambda,\nu)}\nu^+$. Recall
that the beads on the $e$--abacus for $\lambda$ are naturally indexed
by the rows of $\lambda$ and that the $e$--residue of a bead is
defined to be the $e$--residue of the node which is at the end of the
corresponding row. Therefore, an addable $i$--node of $\lambda$
corresponds to a node on runner~$r-1$ of the $e$--abacus which can be
moved to the adjacent position on runner $r$ (which must therefore be
empty).  Similarly, a removable $i$--node corresponds to a node on
runner $r$ which can be moved back to the adjacent position on runner
$r-1$; here, $r\pm 1$ is to be understood modulo~$e$. The addable and
removable nodes of $\lambda^+$ have analogous descriptions.

Fix a partition $\nu$ with $\lambda\iarrow\nu$ and write
$[\nu]=[\lambda]\cup\{x\}$.  Then there exists a node at position
$\beta_x$ on the runner $r-1$ of the $e$--abacus for $\lambda$ which
can be moved to the adjacent position on runner $r$ so as to give the
$e$--abacus for $\nu$. Then $N_i(\lambda,\nu)=A-B$, where
$A=\#\set{y\in A_i(\lambda)|y\succ x}$ and $B=\#\set{y\in
R_i(\lambda)|y\succ x}$.  If $y$ is an addable or removable node of
$\lambda$ then $y\succ x$ if and only if it corresponds to a bead at
position $\beta_y$ with $\beta_y>\beta_x$. Hence, $A$ is equal to the
number of beads on runner $\rho_{r-1}$ which come after $\beta_x$ such
that the adjacent position on $\rho_r$ is vacant; similarly, $B$ is
equal to the number of beads on $\rho_r$ which are after $\beta_x$
and for which the adjacent position on runner~$r-1$ is vacant. 

Now consider the $(e+1)$--abacuses for $\lambda^+$ and $\nu^+$. Assume
first that $r<\alpha$. Then the runners $\rho_{r-1}$ and $\rho_r$ for
$\lambda$ are the same as the runners $\rho^+_{r-1}$ and $\rho^+_r$
for $\lambda^+$ and so the last paragraph shows that the addable and
removable $i$--nodes for $\lambda$ correspond exactly to the addable
and removable $j$--nodes for $\lambda^+$. Hence,
$N_j(\lambda^+,\nu^+)=A-B=N_i(\lambda,\nu)$. Similarly, when
$r>\alpha$ the addable and removable $i$--nodes for $\lambda$
correspond to the addable and removable $(j+1)$--nodes for $\lambda^+$
and $N_i(\lambda,\nu)=N_{j+1}(\lambda^+,\nu^+)$.

Finally, consider the case when $r=\alpha$. This time runner
$\rho_{r-1}$ is equal to $\rho^+_{r-1}$ and runner $\rho_r$ is equal
to $\rho^+_{r+1}$; whereas runner $\rho^+_r=\rho^+_\alpha$ of
$\lambda^+$ is empty. Therefore, the addable and removable $i$--nodes
of $\lambda$ again correspond to addable and removable $j$--nodes
of $\lambda^+$ except that this time there are additional addable
$j$--nodes of $\lambda^+$ corresponding to the adjacent pairs of beads
on the runners $\rho_{r-1}$ and $\rho_r$ of the $e$--abacus of
$\lambda$.  Let $\sigma$ be the partition such that 
$\lambda^+\myarrow{\ j}\sigma\myarrow{j+1}\nu^+$. Since
$\rho_r^+$ is empty, $\lambda^+$ has no removable $j$--nodes.
Therefore, if we let $l$ be the number of pairs of adjacent beads on
runners $\rho_{r-1}$ and $\rho_r$ which are above $\beta_x$ then
$N_j(\lambda^+,\sigma)=A+l$. Next observe that $\sigma$ has a
single addable $(j+1)$--node (corresponding to the bead which we
just moved), and that the removable $(j+1)$--nodes of~$\sigma$
correspond to the removable $i$--nodes of $\lambda$ together with the
$l$ beads on runner~$\rho_r$ which we have already paired with an
adjacent bead on $\rho_{r-1}$; therefore,
$N_{j+1}(\sigma,\nu^+)=-(B+l)$.  Consequently,
$N_j(\lambda^+,\sigma)+N_{j+1}(\sigma,\nu^+)=A-B=N_i(\lambda,\nu)$ and
so we have $$F^+_{j+1}F^+_j\lambda^+ =\sum_{\nu^+}
v^{N_i(\lambda,\nu)}\nu^+ =\Theta( F_i\lambda ),$$ where the sum is
over those partitions $\nu^+$ for which there exists a partition
$\sigma$ such that $\lambda^+\myarrow{\ j}\sigma\myarrow{j+1}\nu^+$.
Note that there are additional terms in the expansion 
of~$F_j\lambda^+$ (corresponding to the pairs of adjacent beads on
runners $\rho_{r-1}$ and $\rho_r$ of the $e$--abacus for $\lambda$);
however, they all disappear when $F_{j+1}$ is applied because these
extra partitions do not have any addable $(j+1)$--nodes. This
completes the proof.  
\end{proof}

Now consider $L(\Lambda_0)_k=L(\Lambda_0)/\(L(\Lambda_0)\cap\F_k)$.
If $\mu$ is an $e$--regular partition with $\len(\mu)\le k$ let
$\tilde B_\mu=B_\mu+\F_{>k}$. As noted by Goodman and
Wenzl~\cite[Lemma~4.1]{GW}, the elements $\set{B_\mu|\mu \text{\ is
$e$--regular and\ }\len(\mu)\le k}$ give a basis of $L(\Lambda_0)_k$.

The bar involution induces a well--defined map on $\F_k$ via
$\bar{a+\F_{>k}}=\bar a+\F_{>k}$ for all $a\in\F$. It is easy to see
that $\tilde B_\mu$ is the unique element of~$\F_k$ which is bar
invariant and of the form $\mu+\sum_\lambda b_{\lambda\mu}(v)\lambda$
for some polynomials $b_{\lambda\mu}(v)\in v\Z[v]$, the sum being over
the partitions of length at most $k$.  

Similarly, 
$\set{\tilde B^+_\nu|\nu 
           \text{\ is $(e+1)$--regular and $\len(\nu)\le k$}}$,
where $\tilde B^+_\nu= B^+_\nu+\F^+_{>k}$, is a basis of the 
$\Umep$--module
$L(\Lambda_0^+)_k=L(\Lambda_0^+)/\(L(\Lambda_0^+)\cap\F^+_k)$.

\begin{Proposition}
Suppose that $\mu$ is an $e$--regular partition with at most $k$
rows. Then $\tilde B^+_{\mu^+}=\Theta(\tilde B_\mu)$.
\label{theta}
\end{Proposition}

\begin{proof}
Looking at the definitions, $\Theta(\tilde B_\mu)=\mu^+$ plus a
$q\Z[q]$--linear combination of less dominant terms. Therefore, it is
enough to show that $\Theta(\tilde B_\mu)$ is a bar invariant element
of $\F^+_k$. 

Let $\0\in\F_k$ be the image of the empty partition in $\F_k$.
Following Lascoux, Leclerc and Thibon~\cite[Lemma~6.4]{LLT} let
$(r_1^{a_1},\dots,r_s^{a_s})$ be the $e$--residue sequence of~$\mu$
corresponding to the $e$--ladders in the diagram of $\mu$. Then
$A_\mu=F_{r_s}^{(a_s)}\dots F_{r_1}^{(a_1)}\0$ is a bar invariant
element of $\F_k$ of the form $A_\mu=\mu+\sum_\lambda
a_{\lambda\mu}(v)\lambda$ where $a_{\lambda\mu}(v)\in\Z[v,v^{-1}]$ and
the sum is over partitions $\lambda$ of length a most $k$ such that
$\mu\gdom\lambda$. Therefore, there exist uniquely determined polynomials
$\alpha_{\sigma\mu}(v)\in\Z[v]$ such that
$\tilde B_\mu=A_\mu-\sum_\sigma \alpha_{\sigma\mu}(v)\tilde B_\sigma,$
where the sum is over $e$--regular partitions $\sigma$ such that
$\mu\gdom\sigma$ and $\len(\sigma)\le k$.

Now consider the element $A_\mu^+=\Theta(A_\mu)=\mu+\sum_\lambda
a_{\lambda\mu}(v)\lambda^+$ in $\F^+_k$. By the Lemma,
$A_\mu^+=\aF_{r_s}^{(a_s)}\dots\aF_{r_1}^{(a_1)}\0$; hence, $A^+_\mu$
is bar invariant. By induction on dominance 
$\tilde B_\sigma^+=\Theta(\tilde B_\sigma)$ for $\mu\gdom\sigma$. 
Therefore, the element
$\Theta(\tilde B_\mu)
  =A^+_\mu-\sum_\sigma \alpha_{\sigma\mu}(v)\tilde B_{\sigma^+}$
is also bar invariant. Consequently, 
$\tilde B^+_{\mu^+}=\Theta(\tilde B_\mu)$ as we wanted to show.
\end{proof}

It is easy to see ~\cite{LLT} that the polynomials
$b^+_{\sigma\tau}(v)$ are non--zero only if $\sigma$ and $\tau$ have
the same $(e+1)$--core.  Therefore, if $\mu$ is an $e$--regular
partition with $\len(\mu)\le k$ then
$\tilde B^+_{\mu^+}
        =\sum_{\mu\gedom\lambda}b^+_{\lambda^+\mu^+}(v)\lambda^+;$
on the other hand,
$\Theta(B_\mu)=\sum_{\mu\gedom\lambda}b_{\lambda\mu}(v)\lambda^+,$
where in both sums $\len(\lambda)\le k$. Hence, Theorem~\ref{fock}
follows from the Proposition.

\section{The main theorem}In this section we extend
Theorem~\ref{fock} to the case where~$\mu$ is not necessarily
$e$--regular; this will prove Theorem~\ref{main}.

The Fock space also admits an action from a {\sf Heisenberg algebra}
$\Hei$~\cite{LT:Schur2}.  The action of $\Hei$ on $\F$ commutes with
the action of $\Ue$ and it is useful because $\F$ is irreducible when
considered as a module for the algebra generated by the actions
of~$\Ue$ and $\Hei$ on~$\F$. In addition, the action of $\Hei$ allowed
Leclerc and Thibon \cite[\S7.9]{LT:Schur2} to extend the bar
involution to the whole of $\F$; in turn, this enabled them to extend
the canonical basis of $L(\Lambda_0)$ to give a basis
$\set{B_\mu|\mu\text{ a partition}}$ of~$\F$ where, the element
$B_\mu$ is again uniquely determined by the two conditions that 
$\bar B_\mu=B_\mu$ and
$$B_\mu=\sum_{\substack{\lambda\vdash n\\\mu\gedom\lambda}}
                 b_{\lambda\mu}(v)\lambda$$
for some polynomials $b_{\lambda\mu}(v)\in\Z[v]$ such that
$b_{\mu\mu}(v)=1$ and $b_{\lambda\mu}(v)\in v\Z[v]$ whenever
$\lambda\ne\mu$. We will show that Theorem~\ref{fock} generalizes to the
non--regular case.

As in the previous sections we are only interested in the action of a
subalgebra~$\Hm$ of $\Hei$; for the full story see
\cite[\S7.5]{LT:Schur2}. The algebra $\Hm$ is generated by elements
$\V_m$ for $m\ge0$; before we can describe how $\V_m$ acts on $\F$ we
need some more notation.

An {\sf $e$--ribbon}  is a connected strip of $e$--nodes which does not
contain a $2\times2$ square; more precisely, an $e$--ribbon is a set
of $e$ nodes $R=\{(a_1,b_1),\dots,(a_e,b_e)\}$ such that
$(a_{i+1},b_{i+1})$ is either $(a_i+1,b_i)$ or $(a_i,b_i-1)$, for
$i=1,\dots,e-1$. The {\sf head} of $R$ is the node $\head(R)=(a_1,b_1)$
and $\spin_e(R)=\#\set{1\le i<e|a_{i+1}=a_i+1}$ is the {\sf $e$--spin}
of $R$.

If $\lambda$ and $\nu$ are partitions then we write
$\lambda\ribbon{m:e}\nu$ if $[\lambda]\ss[\nu]$ and
$[\nu]\setminus[\lambda]$ is a disjoint union of $m$ $e$--ribbons such
that the head of each ribbon is either in the first row of $\lambda$
or is of the form $(i,j)$ where $(i-1,j)\in[\lambda]$. Lascoux,
Leclerc and Thibon (see \cite[\S4.1]{LT:Schur2}), call $\nu/\lambda$
an $e$--ribbon tableau of weight $(m)$ and they note that there is a
unique way of writing $[\nu]\setminus[\lambda]$ as a disjoint union of
ribbons; we will see this below when we reinterpret ribbons in terms
of abacuses. Finally, if~$\lambda\ribbon{m:e}\nu$ then 
$\spin_e(\nu/\lambda)$, the $e$--spin of $\nu/\lambda$, is the sum of
the $e$--spins of the ribbons in $[\nu]\setminus[\lambda]$.

For example, if $\lambda=(3)$ and $e=2$ then the partitions
$\nu$ with $\lambda\ribbon{2:2}\nu$ are
$$\psset{unit=3mm}
  \psline(0,2)(3,2)(3,3)(0,3)(0,2)
  \psline(1,2)(1,3)\psline(2,2)(2,3)
  \psline[linecolor=gray](3,2)(7,2)(7,3)(3,3)
  \psline[linecolor=gray](5,2)(5,3)
\hspace*{21mm}\qquad\qquad
  \psline(0,2)(3,2)(3,3)(0,3)(0,2)
  \psline(1,2)(1,3)\psline(2,2)(2,3)
  \psline[linecolor=gray](3,2)(5,2)(5,3)(3,3)
  \psline[linecolor=gray](0,2)(0,1)(2,1)(2,2)
\hspace*{15mm}\qquad\qquad
  \psline(0,2)(3,2)(3,3)(0,3)(0,2)
  \psline(1,2)(1,3)\psline(2,2)(2,3)
  \psline[linecolor=gray](0,2)(0,0)(2,0)(2,2)
  \psline[linecolor=gray](1,2)(1,0)
\hspace*{9mm}\vrule height 12mm width 0pt
$$
with spins $0$, $0$ and $2$ respectively.

The algebra $\Hm$ is the subalgebra of $\Hei$ generated by elements
$\V_m$ for $m\ge1$. For each $m$, $\V_m$ acts on the Fock space as the
$\C[v,v^{-1}]$--linear map determined by
$$\V_m\lambda=\sum_{\lambda\sribbon{m:e}\nu}
              (-v)^{-\spin_e(\nu/\lambda)}\nu$$
for all partitions $\lambda$. Observe that $\F_{>k}$ is a
$\Hm$--module; hence, there is a well--defined action of $\V_m$
on the quotient space $\F_k$. 

Similarly, there is an action of the negative Heisenberg algebra
$\Hmp$ on the Fock space $\F^+$ and this induces an action on
$\F^+_k$. We denote the generators of $\Hmp$ by $\V^+_m$
for $m\ge1$.

\begin{Lemma}
\label{commuting}
Suppose that $\lambda$ and $\nu$ are partitions of $n$ of
length at most $k$.
\begin{enumerate}
\item If $\lambda\ribbon{m:e}\nu$ then $\lambda$ and $\mu$ have 
the same $e$--core.
\item We have $\lambda\ribbon{m:e}\nu$ if and only if 
$\lambda^+\ribbon{m:e+1}\mu^+$; moreover, if $\lambda\ribbon{m:e}\nu$
then $\spin_e(\nu/\lambda)=\spin_{e+1}(\nu^+/\lambda^+)$.
\end{enumerate}
\end{Lemma}

\begin{proof}
The Lemma will follow once we reinterpret the condition
$\lambda\ribbon{m:e}\nu$ in terms of abacuses. Suppose that
$\lambda\ribbon{m:e}\nu$. Then $[\nu]\setminus[\lambda]$ is a disjoint
union of $e$--ribbons. Extend the partial order $\succ$ on the set of
nodes to a total order by defining $(a,b)\succ(c,d)$ if either $c>a$
or $c=a$ and $b>d$. Totally order the ribbons $R_1,\dots,R_m$ 
in~$[\nu]\setminus[\lambda]$ so that $i>j$ whenever
$\head(R_i)\succ\head(R_j)$. Then the condition that the head of $R_i$
is of the form $(a,b)$ with either $a=1$ or $(a-1,b)\in[\lambda]$ is
equivalent to saying that $[\nu]\setminus(R_1\cup\dots\cup R_i)$ is
the diagram of a partition for $i=1,\dots,m$. Hence, it is enough to
treat the case $m=1$. So let $R=R_1$ where $[\nu]=[\lambda]\cup R$. 

Now the ribbon $R$ is a rim hook and it is well--known that removing a
rim hook of length $e$ from $\nu$ is the same as moving a bead $\beta$
on an $e$--abacus for $\nu$ to the (empty) bead position on the same
runner which is in the preceding row. Further, by definition,
$\spin_e(\nu/\lambda)$ is the leg length of $R$ minus one and, in
terms of the $e$--abacus, the leg length of $R$ is equal to the number
of beads on the abacus which are between the old and new positions
of~$\beta$.  Similarly, the condition $\lambda^+\ribbon{m:e+1}\mu^+$
depends  only on the $(e+1)$--abacuses of $\nu^+$ and $\lambda^+$. As
the $e$ and $(e+1)$ abacuses differ only by the insertion of an empty
runner, all of the assertions of the Lemma now follow.
\end{proof}

\begin{Corollary}
\label{commuting2}
Suppose that $m\ge1$. Then the following diagram commutes.
$$\begin{array}{c}\\[-6pt]
\psmatrix[colsep=1cm,rowsep=1cm]
  \F_k & \F_k \\
  \F^+_k     & \F^+_k
  \psset{nodesep=3pt,arrows=->}
  \everypsbox{\scriptstyle}
  \ncline{1,1}{1,2}^{\V_m}
  \ncline{2,1}{2,2}_{\V^+_m}
  \ncline{1,2}{2,2}\trput{\Theta}
  \ncline{1,1}{2,1}\tlput{\Theta}
\endpsmatrix\\[5pt]
\end{array}$$
\end{Corollary}

\begin{proof}
As with Lemma~\ref{commuting} it suffices to check the result for a
partition $\lambda$ of length at most $k$.
By the definitions and the previous Lemma,
\begin{align*}
\Theta(\V_m\lambda)
  &=\sum_{\lambda\sribbon{m:e}\nu}(-v)^{-\spin_e(\nu/\lambda)}
        \Theta(\nu)\\
  &=\sum_{\lambda^+\sribbon{m:e+1}\nu^+}(-v)^{-\spin_{e+1}
                    (\nu^+/\lambda^+)}\nu^+\\
  &=\V^+_m\lambda^+.
\end{align*}
Therefore, $\V^+_m\(\Theta(\lambda)\)=\Theta\(\V_m \lambda\)$
and we're done.
\end{proof}

Let $\mu$ be a partition. As in the last section let $\tilde B_\mu$
and $\tilde B^+_\mu$ denote the image of $B_\mu\in\F$ and
$B^+_\mu\in\F^+$, respectively,  in $\F_k$ and $\F^+_k$. Then a basis of
$\F_k$ is given by the $\tilde B^\mu$ as $\mu$ runs over all
partitions of length at most $k$ and similarly for $\F^+_k$.

By Proposition~\ref{theta} we know that $\tilde
B^+_{\mu^+}=\Theta(B_\mu)$ whenever $\mu$ is an $e$--regular partition
with $\len(\mu)\le k$.  We can now drop the requirement that $\mu$
should be $e$--regular.

\begin{Proposition}
Suppose that $\mu$ is a partition with $\len(\mu)\le k$.
Then $\tilde B^+_{\mu^+}=\Theta(\tilde B_\mu)$. 
\end{Proposition}

\begin{proof}
As before, the element $\Theta(\tilde B_\mu)$ is equal to
$\mu^+$ plus a $q\Z[q]$--linear combination of less dominant terms, so
it is enough to show that $\Theta(\tilde B_\mu)$ is invariant under
the bar involution since $\tilde B^+_{\mu^+}$ is the unique bar
invariant element in $\F^+_k$ of this form.

By \cite[Prop.~7.6]{LT:Schur2} the bar involution on $\F$ is
completely determined by the conditions $\bar\0=\0$, 
$\bar{F_i^{(m)}x}=F_i^{(m)}\bar x$ and $\bar{\V_m x}=\V_m\bar x$, 
for all $x\in\F$, $0\le i<e$ and $m\ge1$. For each partition 
$\tau=(\tau_1,\dots,\tau_s)$ let
$\V_\tau=\V_{\tau_1}\dots\V_{\tau_s}$. Then
$$\F=\bigoplus_{\tau}\Ume\V_\tau\0$$ 
is a decomposition of $\F$ into a direct sum of irreducible
$\Ue$--modules (where~$\tau$ runs over all partitions of all integers);
see \cite[\S7.5]{LT:Schur2}. Moreover, the modules
$\Ume\V_\tau\0$, for different $\tau$, are all isomorphic as
$\Ume$--modules. Therefore, there exists a bar invariant
basis of $\F$ of the form $A_{\sigma\tau}=F_\sigma\V_\tau\0$ where
$F_\sigma\in\UmeA$, $\sigma$ is an $e$--regular partition and $\tau$
is an arbitrary partition (the elements~$F_\sigma$ are defined in
terms of $e$--residue sequences as in the proof of Proposition~\ref{theta}).
Consequently, we can write 
$\tilde B_\mu=\sum_{\sigma,\tau}a_{\sigma\tau}(v)A_{\sigma\tau}$ for
some bar invariant Laurent polynomials
$a_{\sigma\tau}(v)\in\Z[v,v^{-1}]$.  Now,
$A^+_{\sigma\tau}=\Theta(A_{\sigma\tau})
                 =\Theta(F_\sigma\V_\tau\0)
                 =\aF_\sigma\V_\tau^+\0$ 
by Lemma~\ref{commuting1} and Corollary~\ref{commuting2}; therefore,
$A^+_{\sigma\tau}$ is a bar invariant element of~$\F^+_k$. Hence,
$\Theta(\tilde B_\mu)
         =\sum_{\sigma,\tau}a_{\sigma\tau}(v)A^+_{\sigma\tau}$
is a bar invariant element of~$\F^+_k$, as we needed to show.
\end{proof}

\begin{Remark}
For each composition $\tau=(\tau_1,\dots,\tau_s)$ Leclerc and
Thibon~\cite{LT:Schur2} show that the action of the element
$\V_\tau=\V_{\tau_1}\dots\V_{\tau_s}$ upon~$\F$ is described by
certain polynomials associated with the ribbon tableaux of weight
$\tau$. This is completely analogous to the way in which the action of
$F_{r_s}^{(a_s)}\dots F_{r_1}^{(a_1)}$ on $\F$ can be described in
terms of polynomials associated with standard tableaux.
\end{Remark}

Comparing the coefficient of $\lambda^+$ in $\tilde B^+_{\mu^+}$ and
$\Theta(\tilde B_\mu)$ we obtain the following generalization of our
main theorem.

\begin{Theorem}
Suppose that $\lambda$ and $\mu$ are partitions of length at most
$k$. Then 
$$b_{\lambda\mu}(v)=b^+_{\lambda^+\mu^+}(v).$$
\label{polys}\end{Theorem}

In order to compute the polynomials $b_{\lambda\mu}(v)$ when $\mu$ is
not $e$--regular it is necessary to first invert the ``$R$--matrix''
which describes the bar involution on the basis of~$\F$ given by the
set of partitions. Computationally, this is quite time consuming; in
comparison the regular case is much easier, being essentially Gaussian
elimination. Corollary~\ref{polys} therefore gives a slightly more
efficient way of computing the polynomials
$b_{\lambda\mu}(v)=b^+_{\lambda^+\mu^+}(v)$ since $\mu^+$ is an
$(e+1)$--regular partition by Lemma~\ref{bookkeeping}(i).

Recall that $W^\lambda_q$ and $L^\mu_q$ are the Weyl modules and
simple modules, respectively, for the $q$--Schur algebra
$\Sc_{\C,q}(n,n)$.  Varagnolo and Vasserot~\cite{VV} have shown that
$[W^\lambda_q:L^\mu_q]=b_{\lambda\mu}(1)$; similarly,
$[W^{\lambda^+}_{q'}:L^{\mu^+}_{q'}]=b^+_{\lambda^+\mu^+}(1)$. 

\begin{Corollary}
Suppose that $\lambda$ and $\mu$ are partitions of length at
most~$k$.  Then 
$$[W^\lambda_q:L^\mu_q]=[W^{\lambda^+}_{q'}:L^{\mu^+}_{q'}]
                      =[S^{\lambda^+}_{q'}:D^{\mu^+}_{q'}].$$
\end{Corollary}

\begin{proof}
That
$[W^\lambda_q:L^\mu_q]=[W^{\lambda^+}_{q'}:L^{\mu^+}_{q'}]$ follows
directly from Theorem~\ref{polys} and the remarks above.
For the second claim, observe that the partition $\mu^+$ is
$(e+1)$--regular by Lemma~\ref{bookkeeping}(i); therefore, 
$D^{\mu^+}_{q'}\ne0$. Consequently, 
$[W^{\lambda^+}_{q'}:L^{\mu^+}_{q'}]
                      =[S^{\lambda^+}_{q'}:D^{\mu^+}_{q'}]$
by Schur--Weyl reciprocity.
\end{proof}

Standard Schur functor arguments yield the corresponding statements 
for the $q$--Schur algebras $\Sc_{\C,q}(n,r)$ and $\Sc_{\C,q'}(n,r)$; we leave the
details to the reader.

The last result is interesting because it shows that every
decomposition number for $\Sc_{\C,q}(n,n)$ is also a decomposition
number for some Hecke algebra $\H_{\C,q'}(\Sym_m)$. In contrast
Erdmann~\cite{Erd} has shown that in a given characteristic knowing
all decomposition numbers for the classical Schur algebras
(i.e.~$q=1$) is equivalent to knowing all decomposition numbers for
the symmetric groups (for all $n$ and for a fixed $p$).
Leclerc~\cite{Lec:DecCan} has proved the analogous result relating the
decomposition numbers of the $q$--Schur algebras $\Sc_{\C,q}(n,n)$ and
the Hecke algebras $\H_{\C,q}(\Sym_n)$ (for all $n$ and for a fixed
$q$). No such result is known in the cross characteristic case (i.e.
positive characteristic with $q\ne1$).

Finally, we remark that the full action of $\Ue$ on $\F_k$ and $\Uep$
on $\F^+_k$ are compatible via the map $\Theta$ (in order to make the
statement for $\Upe$ precise~$\F_k$ must be considered as a submodule
of $\F$, rather than a quotient). This can be proved using similar
arguments or, more simply, by invoking \cite[Prop.~7.9]{LT:Schur2}
which says that the actions of $\Ume$ and $\Upe$ on~$\F$ are adjoint
with respect to a natural scalar product on $\F$. The same argument
also proves the corresponding statements for the Heisenberg
algebras~$\Hei$ and $\mathcal H_{e+1}$.

\section{Examples}

Below we give part of the ``crystallised'' decomposition matrices
$\(b_{\lambda\mu}(v)\)$ of the $q$--Schur algebras $\Sc_{\C,q}(n)$
for $(e,n)=(2,6), (3,11), (4,16)$ and $(5,21)$.  By our results,
taking $k=4$ and $\alpha=2$, these submatrices are all the same. 
$$\begin{array}{*5l*{9}{c}}
\cdots&e=5     &e=4    &e=3    &e=2\\\cmidrule[1pt]{1-5}
      &18,3    &14,2   &10,1   &6      &1    &&&&&&&&\\
      &17,4    &13,3   &9,2    &5,1    &v    &1&&&&&&&\\
\cdots&13,8    &11,5   &7,4    &4,2    &.    &v    &1&  &&&&&\\
      &13,4^2  &10,3^2 &7,2^2  &4,1^2  &v    &v^{2}&v&1  &&&&&\\
      &12,9    &9,7    &6,5    &3^2    &.    &.    &v&. &1&&&&\\
      &12,4^2,1&9,3^2,1&6,2^2,1&3,1^3  &v^{2}&v    &v^{2}&v&v&1&&&\\
\cdots&8^2,5   &6^2,4  &4^2,3  &2^3    &.    &.    &v^{2}&v&v&.&1&&\\
      &8^2,4,1 &6^2,3,1&4^2,2,1&2^2,1^2&.    &v^{2}&v^{3}&v^{2}&v^{2}&v&v&1\\
\end{array}$$
Setting $v=1$ we recover the decomposition matrices of the
corresponding $q$--Schur algebras.  Note that when $e>2$ all of the
rows are indexed by partitions which are $e$--regular; therefore, in
these cases the matrix above is a submatrix of the decomposition
matrix for the corresponding Iwahori--Hecke algebra
$\H_{\C,q}(\Sym_n)$; in particular, setting $v=1$ we recover one of
the decomposition matrices from the introduction.
 
To emphasize the dependence on $k$ we again start with $(e,n)=(2,6)$
but now take $k=6$ (and $\alpha=2$); this yields the following
matrices.  {\small $$\begin{array}{*4l*{11}{c}} e=5     &e=4    &e=3
&e=2\\\cmidrule[1pt]{1-4}\\[-11pt] 21,6,3^2&   16,4,2^2&   11,2,1^2&
6      &1&&&&&&&&&&\\ 20,7,3^2&   15,5,2^2&   10,3,1^2&  5,1
&v&1&&&&&&&&&\\ 16,11,3^2&  12,8,2^2&    8,5,1^2&  4,2
&.&v&1&&&&&&&&\\ 16,7^2,3&   12,5^2,2&    8,3^2,1&  4,1^2
&v&v^{2}&v&1&&&&&&&\\ 15,12,3^2&  11,9,2^2&    7,6,1^2&  3,3
&.&.&v&.&1&&&&&&\\ 15,7^2,4&   11,5^2,3&    7,3^2,2&  3,1^3
&v^{2}&v&v^{2}&v&v&1&&&&&\\ 11^2,8,3&    8^2,6,2&    5^2,4,1&  2^3
&.&.&v^{2}&v&v&.&1&&&&\\ 11^2,7,4&
8^2,5,3&5^2,3,2&2^2,1^2&.&v^{2}&v^{3}&v^{2}&v^{2}&v&v&1&&&\\
11,7^2,4^2&8,5^2,3^2&
5,3^2,2^2&2,1^4&v^{2}&v^{3}&.&v&.&v^{2}&.&v&1&&\\
10,7^2,4^2,1&7,5^2,3^2,1&4,3^2,2^2,1&
1^6&v^{3}&.&.&v^{2}&.&.&.&.&v&1\\ \end{array}$$ }%

A consequence of Lemma~\ref{bookkeeping} is that all of these matrices
are the rows of decomposition matrices of the corresponding blocks
which are indexed by partitions with at most $k$ rows (where the
partitions indexing the rows are ordered in a way compatible with
dominance).


\section*{Acknowledgement}
We would like to thank Fred Goodman for some useful discussions.


\let\em\it


\end{document}